\documentclass{article}

\usepackage{PRIMEarxiv}

\usepackage[utf8]{inputenc} 
\usepackage[T1]{fontenc}    
\usepackage{hyperref}       
\usepackage{url}            
\usepackage{booktabs}       
\usepackage{amsfonts}       
\usepackage{nicefrac}       
\usepackage{microtype}      
\usepackage{lipsum}
\usepackage{fancyhdr}       
\usepackage{graphicx}       
\graphicspath{{media/}}     

\usepackage{amsmath,amssymb,amsfonts}
\usepackage{algorithmic}
\usepackage{graphicx}

\DeclareMathOperator*{\argmin}{arg\,min}
\usepackage{algorithm}
\usepackage{xcolor}
\usepackage[english]{babel}

\usepackage{comment}

\usepackage{tkz-euclide}
\usepackage{pst-plot}

\usepackage{caption}
\usepackage{subcaption}
\usepackage{multirow}

\usepackage{tikz}
\usepackage{pgfplots}
\pgfplotsset{compat = newest}

\usetikzlibrary{spy}

\title{GAN-based iterative motion estimation in HASTE MRI}

\author{
  Mathias S. Feinler, Bernadette N. Hahn \\
  Department of Mathematics \\
  University of Stuttgart \\
  Germany\\
  \texttt{\{mathias.feinler, bernadette.hahn\}@imng.uni-stuttgart.de} \\
  }

\begin{document}
\maketitle

\begin{abstract}
	Magnetic Resonance Imaging allows high resolution data acquisition with the downside of motion sensitivity due to relatively long acquisition times. Even during the acquisition of a single 2D slice, motion can severely corrupt the image. Retrospective motion correction strategies do not interfere during acquisition time but operate on the motion affected data. Known methods suited to this scenario are compressed sensing (CS), generative adversarial networks (GANs), and explicit motion estimation. In this paper we propose an iterative approach which uses GAN predictions for motion estimation. The motion estimates allow to provide data consistent reconstructions and can improve reconstruction quality and reliability. With this approach, a clinical application of motion estimation is feasible without any further requirements on the acquisition trajectory i.e. no temporal redundancy is needed. We evaluate our proposed supervised network on motion corrupted HASTE acquisitions of brain and abdomen.

\end{abstract}


\section{Introduction}
Suitable algorithms for motion correction are highly relevant for medical imaging applications such as Magnetic Resonance Imaging (MRI). During data acquisition, patients are asked to remain still. Nevertheless, involuntary motion can occur due to non-cooperative patients or induced by respiratory and cardiac motion. Usually, this comprises motion artefacts in the reconstructions and might lead to an erroneous diagnosis.

Retrospective motion estimation algorithms that are able to compensate motion artefacts with high reliability and accuracy usually require certain redundancies during measurements\cite{cordero2016,burger2018,feinler2023}. Unfortunately, such sequences are not standard in clinical MRI practice. Despite its obvious advantages, acquisition and reconstruction times are longer than for standard sampling schemes. Hence, most radiologists decide to use well-tried cartesian sequences and determine by visual inspection if reconstructions are degraded by motion artefacts\cite{andre2015}. In the latter case the sequence has to be repeated.

The acquisition efficiency has considerably improved due to the introduction of multiple receiver coils\cite{pruessmann1999}. During a pre-scan, sensitivity maps can be computed\cite{allison2012,uecker2014,ma2015}. These allow sub-sampled sequences during the actual imaging procedure. 

One example of a sub-sampled cartesian trajectory is the Half-Fourier Acquisition Single-shot Turbo spin Echo imaging (HASTE) sequence\cite{patel1997}. This sequence uses the k-space symmetry and coil information by an interleaved sampling such that no expressive redundancies remain. The static reconstruction can therefore (almost) exactly replicate the data. In that sense the space of solutions minimizing the data consistency term is non-degenerate and consists mostly of images with motion artefacts in case of motion corrupted data.

In the latter scenario, Generative Adversarial Networks (GANs) can be used to correct for motion artefacts. These networks use the static reconstruction as input and translate it to a motion artefact free image\cite{usman2020,armanious2018b}. This network is also called the Generator. The crucial idea is to use a second network, called Discriminator, that is trained to decide if an image is generated by that Generator or not. Therefore, essentially, the loss of the Generator is learned. This allows to produce visually appealing results. Unfortunately the tweak in this idea is that if the Generator manages to produce just one single artefact free image by completely disregarding the input, the setup fulfils all requirements to perfection. This phenomenon, called mode collapse\cite{goodfellow2016}, can be relaxed by using conditional GANs. These additionally require that the reconstruction is close to the ground truth. Since the artefact free prediction is not unique, this penalty term can only be used moderately to prevent mode collapse. 
Thus, unfortunately, small features can still vanish or emerge for no apparent reason. Even though the relevant literature\cite{usman2020, yoshida2022, armanious2018b, bao2022, johnson2019} manages to improve image quality, a closer visual inspection reveals that especially the recovered individual brain structures is lacking detail. Bearing in mind that mostly small features are the core to do MRIs in the first place, GANs are not reliable in clinical situations where the detection of small subtle malignancies is crucial. 

If the deformation fields responsible for the corrupted data are known, a simple CG-SENSE plus motion procedure can produce data consistent and motion artefact free reconstructions. Motion estimation therefore is the foundation to reduce motion artefacts while keeping data consistency. In \cite{haskell2019} the authors use a convolutional neuronal network (CNN) to assist rigid motion estimation. The motion parameters are optimized by an analytic Newton procedure. Therefore, the number of iterations is empirically high. Further, the convergence for large deformations is precarious\cite{cordero2016} due to the localized gradient information, and an application to general deformations is infeasible due to the intractable number of deformation parameters. 

Dynamic MRI for the time-resolved measurement of the heart has shown huge interest for more accurate reconstruction methods. If the amount of acquired k-space data per time-frame is significantly reduced to increase temporal resolution, even state-of-the-art sparse reconstruction algorithms\cite{adler2017,adler2018,hammernik2018} are unsuccessful. For compensation, the data of other frames are used additionally\cite{schlemper2017}, bearing in mind, that these are subject to different deformation states. The easiest way to incoprorate such information is a registration of binned data, to enable a reconstruction using the full data\cite{balakrishnan2019, munoz2022}. Inspired by computer-vision, recurrent all pairs field transforms (RAFT) have been used, to improve the registration quality\cite{seegoolam2019, teed2020, hammernik2021, pan2022}. All these methods use a feature extractor and a context encoder applied to the fixed and moving images, respectively, and compute subsequently 4D correlations and estimate motion. Hence these methods require a \emph{similar} image-representation of the sub-sampled data. Only if consecutively the center of k-space is sampled or a radial/spiral trajectory is used, such methods are likely to converge with good reconstruction quality. 

In contrast, for the HASTE sequence, the convergence of aforementioned methods is unlikely. Data is acquired from the centre of k-space to the periphery. Hence, sub-sampled data is in this case a small stripe of the k-space with disjoint image representations. 

To tackle the above, we propose an unrolled iterative Network using GANs to assist motion estimation. Instead of an analytic gradient descent, we formulate motion estimation as a registration task which is solved using Deep CNNs. We evaluate the proposed procedure on simulated motion corrupted HASTE MRIs of the brain and the abdomen.

The remainder of this paper is organized as follows. Section II introduces the theoretical basis including the forward model, k-space trajectory and the formulation of the motion estimation and image reconstruction task.
In Section III we develop our motion estimation strategy. We introduce a simplified forward model which is computationally feasible, adapt GANs to suit an iterative process and introduce the structure of the proposed network. Further, we derive a suitable training loss which is independent on any reference configuration and present how our dynamic training data is synthetically generated. 
The results on simulated brain and abdominal MRIs are shown in Section IV and a conclusion is drawn in Section V.

\section{Theory}

\subsection{Forward Imaging Model} 
A suitable forward operator is required to model the imaging process. In the dynamic case, this operator maps a reference image $s^\mathrm{ref}$ and a time-continuous deformation field $U^\mathrm{ref}$ to measured MRI data $y$. Essentially, we use a continuous version of the multishot motion model of Batchelor et al.\cite{batchelor2006}
\begin{equation}
	\mathcal{A}_{t,c}(U,s) = M_{t} \mathcal{F}[S_c U_{t}[s]], \label{eq:model} 
\end{equation}
where $t \in [0, T]$ is the time index with $T$ being the total acquisition time, $U \in \mathcal{D}_\mathrm{U}$ is the (a priori unknown) time and space-dependent deformation field contained in the space of admissible deformations $\mathcal{D}_\mathrm{U}$, $s \in \mathcal{D}_\mathrm{s}$ is the searched-for image at reference configuration contained in the space of admissible ground truth images $\mathcal{D}_\mathrm{s}$, $\mathcal{F}$ is the Fourier transform and $M_{t}$ is the characteristic function for the k-space portion extracted at time $t$. The $N^\mathrm{coils}$ coil sensitivity maps $\{S_c\}_{c=1}^{N^\mathrm{coils}}$ are assumed to be known by some precomputation procedure\cite{allison2012} and are time-invariant. $U_{t}[s]$ corresponds to the deformed state of the ground truth image at excitation time $t$. 

The space of admissible deformations $\mathcal{D}_\mathrm{U}$ is a generic space that is implicitly defined in context with each application. For brain MRIs $\mathcal{D}_\mathrm{U}$ might be limited to rigid deformations, whilst for abdominal free breathing MRIs, much more complex deformations define $\mathcal{D}_\mathrm{U}$ containing the totality of possible deformations caused by human breathing.

In general, measured data $y_t^c = \mathcal{A}_{t,c}(U^\mathrm{ref},s^\mathrm{ref}) + \mathcal{N}$ are only available with noise. We model this measurement noise $\mathcal{N}$ as additive normal distributed noise. We can split $y$ into contributions per time and coil sensitivity map by $y = \{y_t\}_{t} = \{y_t^c\}_{t,c}$.

We further define the field of view $\Omega^\mathrm{FOV} $ as a two-dimensional rectangle. 

Using this forward imaging model for dynamic MRI, we aim to estimate deformations $U$ and reconstruct the image $s$, given data $y$ and coil sensitivities $S_c$.

\subsection{Sampling Trajectories}
The choice of the sampling trajectory determines the chronology in which data is acquired. The trajectory therefore determines if and up to which accuracy motion estimation is feasible. 

The HASTE sequence can be written as 
\begin{equation}
	M_t(x,y) = \chi_{\delta_{\pi t/T,y}} , \quad t = 2jT/N, \, j \in [0,1,\ldots,N/2], \label{eq:HASTE}
\end{equation}
where $\delta$ is the Kronecker delta and $\chi_{\delta_{\pi t/T,y}}$ is the characteristic function for the set $\delta_{\pi t/T,y}$. The coordinates of the k-space are normalized to $[-\pi,\pi]^2$. 

The k-space lines are acquired consecutively at discrete points in time, starting in the centre of k-space. The phase encoding (PE) direction specifies the line to be acquired. The frequency encoding (FE) direction describes the position across this line. Hence, in FE-direction the full frequency information is measured, while in PE-direction only one particular frequency is encoded. The complete information in FE-direction leads to sharp edge information in FE-direction, while edges in PE-direction are not represented. Taking multiple consecutive lines into account allows a rough identification of edges in PE-direction. In particular, the number of considered lines has to be balanced. Too many lines allow motion to smear out edges, while to few lines might not be enough to encode edges in the first place. Since edges are crucial for motion estimation, already at this point it is clear that motion estimation in PE-direction will be much more challenging than in FE-direction for the HASTE sequence.

\subsection{Reconstruction and Motion Estimation}
Now that the forward imaging model \ref{eq:model} and the data acquisition sequence \ref{eq:HASTE} is defined, we consider the inverse problem of reconstructing the the tuple $(s, U)$ given data $y$ and coil sensitivities $S_c$. Assuming the measurement noise is normally distributed, the appropriate penalty term for data consistency is given by the squared $L_2$ norm leading to 
\begin{equation}
	(s^*, U^*) = \argmin_{U \in \mathcal{D}_\mathrm{U}, s \in \mathcal{D}_\mathrm{s}} \|y - \mathcal{A}(U, s) \|^2, \label{eq:reco}
\end{equation}
c.f.\cite{scherzer2008}. The restriction $s \in \mathcal{D}_\mathrm{s}$ is crucial for the case of the HASTE sequence. The HASTE sequence does not provide significant redundancy, i.e. the unconstrained static reconstruction
\begin{equation}
	s^\mathrm{stat} = \argmin_{s} \|y - \mathcal{A}(I, s) \|^2 \nonumber
\end{equation}
corresponds to a negligibly small residual value even if $U^\mathrm{ref}$ differs severely from the identity $I$. Therefore, reconstruction and motion estimation must essentially be based on enforcing $s^* \in \mathcal{D}_\mathrm{s}$

\section{Methodology}
\subsection{Motion Estimation using Ground Truths} \label{subsec:motest}
Suppose, the ground truth $s^\mathrm{ref}$ is available, and the only remaining task is to compute the deformation fields $U^\mathrm{ref}$. This problem can be written as
\begin{equation}
	U^* = \argmin_{U \in \mathcal{D}_\mathrm{U}} \|y - \mathcal{A}(U, s^\mathrm{ref})\| .
\end{equation}
For any time instance $t$, we have to find a suitable deformation $U_t$ such that $y_t$ and $\mathcal{A}_t(U_t, s^\mathrm{ref})$ are consistent with each other. An unbiased first guess is $U_t^{(0)} = I$. By using a static reconstruction operator $\mathcal{R}^{\mathrm{stat}}$, we can now compute partial reconstructions $s_t^{(0),\mathrm{G}} = \mathcal{R}^{\mathrm{stat}}[\mathcal{A}_t(I, s^\mathrm{ref})]$ and $s_t = \mathcal{R}^{\mathrm{stat}}[y_t]$. 

These partial reconstructions are image representations of a limited frequency component. Empirically, $s_t$ and $s_t^{(0),\mathrm{G}}$ look alike but are subject to a different configuration. The estimation of a deformation field that morphs one into the other is called a registration task and has gained huge interest for classical\cite{horn1981} and learned methods\cite{ilg2016}.

Unfortunately, the problem is that even without measurement noise, in general, $U_t^\mathrm{ref}[s_t^{(0),\mathrm{G}}] \not= s_t$. Since the sampling pattern is fixed, we essentially compare different frequency information due to the known relation of deformation in image and frequency space\cite{vaillant2014}. However, the deformation state is encoded. Hence, a comparison of $s_t^{(0),\mathrm{G}}$ and $s_t$ allows to compute a first \emph{estimate} $U_t^\mathrm{(1)}$ to $U_t^\mathrm{ref}$. Since we compare different frequency information, we cannot expect that first guess to be exact. But, we can now compute an updated version $s_t^{(1),\mathrm{G}} = \mathcal{R}^{\mathrm{stat}}[\mathcal{A}_t(U^{(1)}, s^\mathrm{ref})]$. 

By comparing $s_t^{(k),\mathrm{G}} = \mathcal{R}^{\mathrm{stat}}[\mathcal{A}_t (U^{(k)}, s^\mathrm{ref}) ]$ and $s_t$, only the remaining inexactness in the deformation has to be measured and corrected. 

Therefore, motion estimation has to be performed in an iterative fashion.

Indeed, as some approximation $U^{(k)}$ converges to $U^\mathrm{ref}$, also $s_t^\mathrm{G}$ converges to $s_t$.

\subsection{Generative Adversarial Networks (GANs) for Motion Correction}\label{subsec:gan}
Since the ground truth $s^\mathrm{ref}$ is unknown, we at least have to find a useful approximation to apply motion estimation. For blind retrospective motion correction, Generative Adversarial Networks (GANs) are the most promising approach due to fast computation times and universal applicability, even though data consistency is not guaranteed. GANs consisting of a Generator (G) and a Discriminator (D) are usually trained in a min-max game with the loss 
\begin{equation}
	\min_G \max_D \,\, \mathrm{log}(1-D(G(s^\mathrm{stat}))) + \mathrm{log}(D(s^\mathrm{ref})) \nonumber \label{eq:LossGAN}
\end{equation}
For the purpose of motion correction, conditional GANs are used and trained with the additional loss 
\begin{equation}
	\psi \|s^\mathrm{ref} - G(s^\mathrm{stat})\|_2^2,\label{eq:gan-ref-stat}
\end{equation}
where $\psi \in [0, \infty)$ balances the contributions\cite{usman2020}. Due to non-uniqueness, the training with $\psi \rightarrow \infty$ does not produce sharp results, while $\psi=0$ intensifies mode collapse and therefore unreliable predictions.

For the setup of motion correction, we want severe motion artefacts to be corrected well, while unnoticeable motion artefacts are dispensable. Hence, we define the function $\kappa(U) = (1-\|U\|/U_\mathrm{max})^{\frac{1}{5}}$ with $U_\mathrm{max}$ being the maximum over all $\|U\|$ of the considered data set. Using this function we propose 
\begin{equation}
	\min_G \max_D \,\, \kappa(U)\cdot \mathrm{BCE}(1,D(G(s^\mathrm{stat}(U)))) + \mathrm{BCE}(0, D(s^\mathrm{ref})) \label{eq:LossGAN_new}
\end{equation}
where $\mathrm{BCE}$ is the binary cross entropy. This form already includes that static reconstructions without motion are mapped to artefact free images. To achieve a stable steady state at the searched-for solution, we want the Generator to do \emph{nothing} if $s^\mathrm{ref}$ is plugged in. We therefore explicitly add the following contribution to the min-max loss
\begin{equation}
	\psi \|s^\mathrm{ref} - G(s^\mathrm{ref})\|_2^2. \label{eq:gan-ref-ref}
\end{equation}
For empirical reasons, we weight the contributions \eqref{eq:gan-ref-stat} and \eqref{eq:gan-ref-ref} with the same parameter $\psi$. An individual choice is possible. The contributions \eqref{eq:LossGAN_new}, \eqref{eq:gan-ref-stat} and \eqref{eq:gan-ref-ref} define the loss function for our purpose.

\subsection{A Computationally Feasible Forward Imaging Model}\label{subsec:subsection_deformation_model}

The considered deformations $U:\mathbb{R}^2\times[0,T] \rightarrow \mathbb{R}^2$ are vector fields in space and time. 
If we deform an image $s:\mathbb{R}^2 \rightarrow \mathbb{C}$ using a vector field, we evaluate the image at the position the vector field points to. Hence we write
\begin{equation}
	U_t[s]\left(\begin{pmatrix}x\\y\end{pmatrix}\right) := s\left(U_t\left(\begin{pmatrix}x\\y\end{pmatrix}\right)\right),
\end{equation}
where time $t$ is indexed as a subscript. The notation consistent with the discretized setting, where deformation is a matrix vector multiplication. 

In the general setting, the parameters of a deformation are the vector field itself. Therefore, at each position and time, a 2D vector needs to be specified. 

To make motion estimation feasible, we have to reduce the temporal complexity. We use a fixed basis of $N^\mathrm{U}$ deformation fields $\{U^H_i\}_{i=1}^{N^\mathrm{U}}$ to define deformations for all times via piece-wise linear interpolation over time, i.e.
\begin{align}
	U_t = &(1-d)U^H_{i} + d U^H_{i+1}, \nonumber \\
	&h = \frac{t}{T} N^\mathrm{U}-\frac{1}{2}, \, i = \lfloor  \min(\max(h,0),N^\mathrm{U}-2) \rfloor , \, d = h-i. \nonumber
\end{align}
The multiplication by a constant is to be understood as a scaling of the respective parametrizations, where the identity deformation is retrieved if all parameters are zero. We want to point out, that deformations are non-linear with respect to the parameters which means $(U_1 + U_2)[s] \not = U_1[s] + U_2[s]$ in general. In an abuse of notation, we used the sum of deformations to indicate the sum of the respective parametrizations.

The forward operator \ref{eq:model} models the imaging process accurately. For reconstruction, however, the computational complexity has to be balanced with accuracy. We therefore introduce an approximative model for reconstructions. 

First, we evaluate $U_t$ at $\eta N^\mathrm{U}$ locations, where $\eta>=1$ is a natural number. This provides us with the deformation fields $\{U^K_j\}_j$. Further, we approximate 
\begin{equation}
	U_t[s] \approx (1-d)U^K_j[s] + dU^K_{j+1}[s],
\end{equation}
where $h = \frac{t}{T} \eta N^\mathrm{U}-\frac{1}{2}, \, j = \lfloor \min(\max(h,0),\eta N^\mathrm{U}-2) \rfloor , \, d = h-j$. Note, that here we linearly combine \emph{evaluated} deformations. This approximation is only for small changes in deformations a suitable choice.
The model for reconstruction reads now
\begin{equation}
	\tilde{y}_t^c = M_{t} ((1-d) \mathcal{F}[S_c U^K_j [s]] + d\mathcal{F}[S_c U^K_{j+1} [s]]).
\end{equation}
Hence, only $\eta N^\mathrm{U}$ deformations and Fourier transforms need to be calculated. The hyper-parameter $\eta$ therefore allows to balance accuracy and computational complexity. 

We want to point out that a piecewise constant ansatz does not work properly for general sequences, since discontinuous changes of deformation fields in time can cause streaking artefacts for reconstructions. Using aforementioned approximations we were able to diminish such artefacts almost entirely.

\subsection{Deformations in Discrete Settings}

In discrete settings, a deformation by non-integer shifts is not uniquely defined. Therefore, interpolation strategies are usually used. In the literature, several such strategies have been compared\cite{unser1995}. The issue with bilinear interpolation becomes clear for the following simple example. If we shift an image by one half pixel to the right, and afterwards one half pixel to the left, we would expect to get back the image we originally had in hand. Unfortunately, the result will be a smoothed version of the original image. 

To prevent such smoothing, higher order interpolation is used. In special cases like rigid transforms, exact theory exists\cite{unser1995}, by using the Fourier shift theorem. For general deformation fields, spline interpolation is mostly used.

\begin{figure}
	\centering
	\begin{subfigure}[b]{0.25\textwidth}
		\centering
		\includegraphics[width=\textwidth]{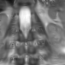}
		\caption*{bilinear}
	\end{subfigure}
	\begin{subfigure}[b]{0.25\textwidth}
		\centering
		\includegraphics[width=\textwidth]{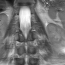}
		\caption*{ground truth}
	\end{subfigure}
	\begin{subfigure}[b]{0.25\textwidth}
		\centering
		\includegraphics[width=\textwidth]{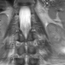}
		\caption*{sinc interpolation}
	\end{subfigure}
	
	\caption{Results of $U^{-1}[U[s]]$ by using different versions for discretized transforms. The bilinear interpolation smooths out details. The sinc interpolation with 7 times 7 patches is visually almost artefact free.}
	\label{fig:general_sinc_trafo}
\end{figure}

Since the Fourier shift theorem provides excellent results for rigid transforms, we study the idea that every non-linear transform can locally be approximated by a linear transform. Therefore we use a small patch of e.g. 7 times 7 pixels and apply the Fourier shift by the remaining non-integer modulus. This is essentially a sinc interpolation. The results show superior quality compared to several other interpolation methods when used with non-linear deformations. We therefore use sinc interpolation for an accurate forward simulation of the dynamics but stick with the bilinear operator for reconstruction due to time-advantage. A small comparison is shown in Figure \ref{fig:general_sinc_trafo}, where we use some deformation and its inverse, to show the effects of $U^{-1}[U[s]]$.

\subsection{Motion estimation using GANs}

In this subsection we introduce an iterative network structure using pre-trained GANs to estimate motion. 

\begin{figure}
	\centering
	\includegraphics[width=0.6\textwidth]{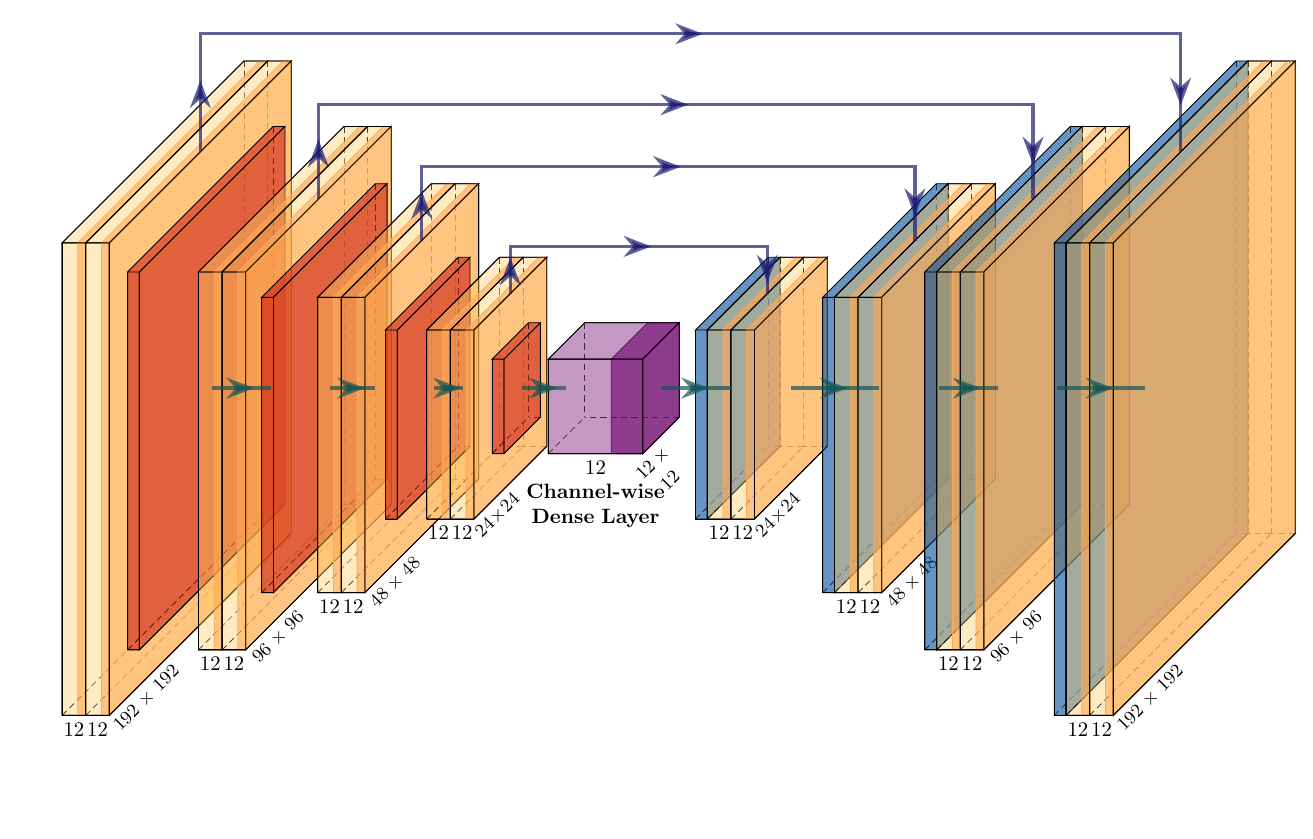}
	\caption{$\text{V-netD}^{^{D,C}}_\Theta$ with parameters $\Theta$, depth $D=4$ and $C=12$ channels for images of resolution $N\times N=192 \times 192$.}
	\label{fig:V-netD}
\end{figure}

The basic network module is inspired by the networks used for inpainting tasks. These are essentially U-net like structures with an additional channel-wise dense layer for better extrapolation and information distribution across the image domain\cite{yang2017}. We refer to this module as V-NetD in the following, depicted in figure \ref{fig:V-netD}. We want to combine multiple V-NetD to a much more complex network structure that incorporates the intuition gained from sections \ref{subsec:motest} and \ref{subsec:gan} on the task of motion estimation. 

\begin{algorithm}
	\caption{iterative GAN}
	\label{alg:iterGAN}
	\begin{algorithmic}[1]
		\STATE{$\mathrm{Input}:  \, \, y_i, \,\, i = 1,\ldots,N^\mathrm{U}$} 
		\STATE{$\mathrm{Output}:  \, \, U_i^\mathrm{est}:= U_i^{(N^{\mathrm{iter}}+1)}, \,\, i=1,\ldots,N^\mathrm{U}$}
		\STATE{$\mathrm{Initialize}: s^{(1)} = 0, \,\, U^{(1)}=I $}
		
		\FOR{ $k=1,\ldots,N^\mathrm{iter}$}
		\STATE{$ s^{(k+1)} = \mathrm{CG}(y,U^{(k)},s^{(k)}, N^\mathrm{CG})$}            
		
		\STATE{$ s^\mathrm{GAN} = \mathrm{GAN}(s^{(k+1)})$}
		
		\STATE{$ y^\mathrm{GAN} = \mathcal{A}(s^\mathrm{GAN}, U^{(k)})$}

		\FOR{ $i = 1, \ldots, N^\mathrm{U}$}   
		
		\STATE{$ s_i = \mathcal{R}^{\mathrm{stat}}[y_i] $}
		\STATE{$ s_i^\mathrm{GAN} = \mathcal{R}^{\mathrm{stat}}[y_i^\mathrm{GAN}]$}
		
		\STATE {$U^{\Delta}_i = {\textbf{\text{V-netD}}^{4,32}_{\mathcal{E}^{(k,i)}}}(s_i, s_i^\mathrm{GAN}, U_i^{(k)}[s^{(k)}], U_i^{(k)}[s^\mathrm{GAN}])$} 
		
		\ENDFOR
		\STATE {$U^{(k+1)} = $ \parbox[t]{0.7\linewidth}{%
				$U^{(k)} + U^{\Delta} +  \mathbin{\textbf{\text{V-netD}}^{4,16N^\mathrm{U}}_{\Psi^{(k)}}}(U^{(k)}, U^{\Delta})$\;}} 
		\ENDFOR
	\end{algorithmic}
\end{algorithm}
The procedure is presented in Algorithm \ref{alg:iterGAN}.
Herein $\mathrm{CG}(y,U,s,N^\mathrm{CG})$ is the CG-SENSE algorithm applied to the dynamic forward operator with data $y$, deformations $U$, initial value $s$ and iteration number $N^\mathrm{CG}$\cite{pruessmann2001}.
The parameter $N^\mathrm{CG}$ can be adapted to the characteristics of the data. In our case we fix the parameters $N^\mathrm{iter}=3$ and $N^\mathrm{CG}=5$. The data $\{y_t\}_t$ is naturally reorganized into $N^\mathrm{U}$ fractions $\{y_i\}_{i=1}^{N^\mathrm{U}}$. 

The network on images of resolution $192\times192$ with $N^\mathrm{U}=8$ has 27.8 Million parameters of which 8.6 Million parameters are fixed due to pre-trained Generators. 

\subsubsection{Network Training of iterative GAN}


For the training of neuronal networks a loss function has to be defined. Most networks are trained in supervised fashion. The dataset requires therefore the exact knowledge of the outcomes. The problem with motion estimation is that reference configurations are in general not unique. For residual based motion estimation, we can fix one data component and align all further measurements by motion estimation. With this approach, the map is one-to-one. In case of no redundancies, the prediction of the GAN can possibly resemble a different reference configuration. Recall, that the requirement of GANs was to produce images that are indistinguishable from originals, but not necessarily exactly equal. Therefore, appropriate deformation fields can be subject to a different reference configuration.

Let $s^\mathrm{ref}$ be the ground truth reconstruction in reference configuration. The according deformation field is given by $U^\mathrm{ref} = \{U_i^\mathrm{ref}\}_i$.
Let the GAN prediction $s^\mathrm{GAN} = U^\mathrm{R}[s^\mathrm{ref}]$ be exact, modulo some deformation $U^\mathrm{R}$. Data consistency is therefore reached when we combine $s^\mathrm{GAN}$ with $\{U_i^\mathrm{ref}[U^\mathrm{R}[\cdot]]\}_i$. Therefore, we have to find some reference configuration independent loss, to train such networks for motion estimation. Using the mean deformation $U^\mathrm{ref,mean} = \frac{1}{N} \sum_{i} U_i^\mathrm{ref}$ we can define 
\begin{align}
	& U_i^\mathrm{ref}[U^\mathrm{R}[ (U^\mathrm{ref,mean}[U^\mathrm{R}[ \cdot ]])^{-1} ]] \\
	=& U_i^\mathrm{ref}[U^\mathrm{R}[ ((U^\mathrm{R})^{-1}[(U^\mathrm{ref,mean})^{-1}[ \cdot ]]) ]] \\
	=& U_i^\mathrm{ref}[(U^\mathrm{ref,mean})^{-1}[ \cdot ] ] =: \tilde{U}_i^\mathrm{ref}
\end{align}
which holds if $U^\mathrm{R}$ is invertible. The crucial property is that this expression is $U^\mathrm{R}$ independent, with the price that we have to compute an inverse  deformation field. Thus, we use the following loss function
\begin{equation}
	\mathcal{L}(U^\mathrm{est}, U^\mathrm{ref}) = \|\tilde{U}^\mathrm{est} - \tilde{U}^\mathrm{ref}\|_2^2 + \beta \|\Delta U^\mathrm{est}\|_2^2, 
\end{equation}
where $\beta>0$ promotes smooth motion estimates.

Recall, that the choice of $U^\mathrm{mean}$ for the computation of a reference configuration independent loss is crucial for several reasons. First, it might seem acceptable to use e.g. $U_1^\mathrm{ref}$ instead. The problem with this choice is that $\nabla_{U_1}\mathcal{L} = 0$. Thus the network training would leave the parameters affecting exclusively $U_1$ unchanged. Second, some symmetry is enforced, since all gradients $\nabla_{U_i}\mathcal{L}$ are given in terms of the same expression. Third, the gradient information is consistent, i.e. only $\nabla_{U_i}\mathcal{J}$ is required for an optimal update of $U_i$.

The computation of the inverse of $U^\mathrm{mean}$ is non-trivial and can be realized by minimizing the functional
\begin{equation}
	\min_{V}\| V[U^\mathrm{mean}[\cdot]] - I[\cdot] \|^2 + \alpha \|\nabla V\|^2,
\end{equation}
where $I$ is the identity operator and $V$ will be the approximation to $(U^\mathrm{mean})^{-1}$. The regularization term, controlled by the parameter $\alpha$, requires smoothness of $V$. A Landweber iteration can approximate the inverse sufficiently well within a few steps.

For the training of the network we employ a cascade training procedure\cite{marquez2018}. Therefore, in contrast to end-to-end training, no gradients across CG steps need to be computed. We start at the first cascade ($k=1$) and work our way up to the last cascade. With this approach, intermediate predictions are only dependent on network parameters of deeper layers. Since these are fixed, also intermediate predictions do not change. Hence, we can precompute the inputs to the cascade under training. With this approach the training is in reasonable time feasible. During training of each cascade, we keep the parameters of the Generator fixed. Thus, motion estimation can be trained in standard supervised form. 

For each cascade, we can retrain the Generator to suit remaining motion artefacts. For all cascades $k>2$ it is crucial, that we initialize the weights with the converged weights of the previous cascade for training. Empirically, the performance improved by approximately 20\%. The reason is that ab initio training on noisy data with small remaining inaccuracies can lead to networks which rather judge if $U^\mathrm{est} \in \mathcal{D}_\mathrm{U}$. This will lead to overfitting. In contrast the initialization already learned to mainly focus on inconsistencies between $s_i$ and $s_i^\mathrm{GAN}$, which extrapolates well to unseen data.

\subsection{Training Data Generation}

\subsubsection{Brain MRI}\label{sec:def_brain_MRI}
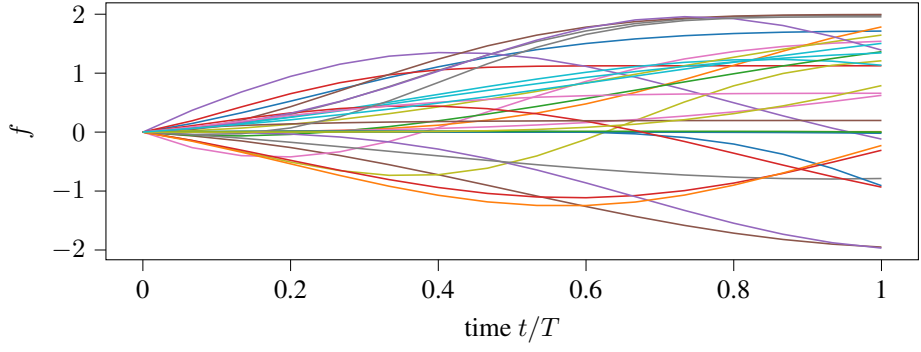
\begin{figure}
	\centering
	\begin{tikzpicture}
		
		\definecolor{crimson2143940}{RGB}{214,39,40}
		\definecolor{darkgray176}{RGB}{176,176,176}
		\definecolor{darkorange25512714}{RGB}{255,127,14}
		\definecolor{darkturquoise23190207}{RGB}{23,190,207}
		\definecolor{forestgreen4416044}{RGB}{44,160,44}
		\definecolor{goldenrod18818934}{RGB}{188,189,34}
		\definecolor{gray127}{RGB}{127,127,127}
		\definecolor{mediumpurple148103189}{RGB}{148,103,189}
		\definecolor{orchid227119194}{RGB}{227,119,194}
		\definecolor{sienna1408675}{RGB}{140,86,75}
		\definecolor{steelblue31119180}{RGB}{31,119,180}
		
		\begin{axis}[
			tick align=outside,
			tick pos=left,
			height=5cm,width = 0.75\textwidth, 
			xmin=-0.075, xmax=1.575, xlabel=time $t/T$,
			xtick style={color=black}, xtick = {0, 1.5/5, 2*1.5/5, 3*1.5/5, 4*1.5/5, 1.5}, xticklabels={0, 0.2, 0.4, 0.6, 0.8, 1},
			y grid style={darkgray176},
			ymin=-2.1658155773519, ymax=2.19326363568076, ylabel=$f$,
			ytick style={color=black}
			]
			\addplot [semithick, steelblue31119180]
			table {%
				0 0
				0.1 0.146608017278622
				0.2 0.326891932443109
				0.3 0.525013799950756
				0.4 0.728579643284352
				0.5 0.92609765101296
				0.6 1.10788397119633
				0.7 1.26686716929839
				0.8 1.39900365640428
				0.9 1.50324600633326
				1 1.58111347195774
				1.1 1.63598559973708
				1.2 1.67227668920768
				1.3 1.69464845865743
				1.4 1.70738498687436
				1.5 1.7139992683259
			};
			\addplot [semithick, darkorange25512714]
			table {%
				0 0
				0.1 0.00191457304246834
				0.2 0.00385870199789362
				0.3 0.00528074819733926
				0.4 0.00616623163579222
				0.5 0.0066400736609139
				0.6 0.00685324408641841
				0.7 0.00693003023238636
				0.8 0.00695025521938275
				0.9 0.00695351555735391
				1 0.00695373033681113
				1.1 0.00695373200710492
				1.2 0.00695372956923035
				1.3 0.00695339236930415
				1.4 0.00694745072959557
				1.5 0.00690353569988311
			};
			\addplot [semithick, forestgreen4416044]
			table {%
				0 0
				0.1 7.69661265314081e-05
				0.2 8.152223663771e-05
				0.3 8.1767014228018e-05
				0.4 8.17724467836456e-05
				0.5 8.17724619532889e-05
				0.6 8.1772461918983e-05
				0.7 8.17723027585204e-05
				0.8 8.17599325418916e-05
				0.9 8.15366440997112e-05
				1 7.96514924650271e-05
				1.1 6.96384339974321e-05
				1.2 3.06751151581963e-05
				1.3 -9.05992070894257e-05
				1.4 -0.000409597000862627
				1.5 -0.00114511490174973
			};
			\addplot [semithick, crimson2143940]
			table {%
				0 0
				0.1 0.191851411307141
				0.2 0.427363119288496
				0.3 0.652061939948188
				0.4 0.838119358871162
				0.5 0.972362564358207
				0.6 1.05573270644397
				0.7 1.09927257033019
				0.8 1.11770006710694
				0.9 1.12366081478836
				1 1.12499735080123
				1.1 1.12517128572593
				1.2 1.12518035851486
				1.3 1.12518042869346
				1.4 1.12518042836628
				1.5 1.12517957718831
			};
			\addplot [semithick, mediumpurple148103189]
			table {%
				0 0
				0.1 0.370016697809627
				0.2 0.680456712343086
				0.3 0.94576818494017
				0.4 1.15198715274569
				0.5 1.28823012819566
				0.6 1.34925383641248
				0.7 1.3358339504256
				0.8 1.25424671979776
				0.9 1.1152149976781
				1 0.932540027557726
				1.1 0.721604005454766
				1.2 0.497904717117179
				1.3 0.275755854388438
				1.4 0.0672532057325802
				1.5 -0.118430129202576
			};
			\addplot [semithick, sienna1408675]
			table {%
				0 0
				0.1 -0.0640879132904466
				0.2 -0.150113965119053
				0.3 -0.261869152312979
				0.4 -0.397346669784315
				0.5 -0.552772839871592
				0.6 -0.723105402763423
				0.7 -0.902440250392112
				0.8 -1.08440339646154
				0.9 -1.2625233767981
				1 -1.43057067460503
				1.1 -1.58285192912205
				1.2 -1.71444984771358
				1.3 -1.82140329445873
				1.4 -1.90082547663141
				1.5 -1.9509612778723
			};
			\addplot [semithick, orchid227119194]
			table {%
				0 0
				0.1 -0.263283026406546
				0.2 -0.399529343464454
				0.3 -0.420433388480154
				0.4 -0.334144318793681
				0.5 -0.159538500280868
				0.6 0.0763856744263252
				0.7 0.343298135167702
				0.8 0.612487961459813
				0.9 0.860785480734427
				1 1.07279940293373
				1.1 1.24133512848811
				1.2 1.3662984222077
				1.3 1.45265946333112
				1.4 1.50812409667129
				1.5 1.54105778869472
			};
			\addplot [semithick, gray127]
			table {%
				0 0
				0.1 0.0423638389932388
				0.2 0.140276533634382
				0.3 0.299511483832254
				0.4 0.515356069037352
				0.5 0.77197031732671
				0.6 1.04534194657406
				0.7 1.30845014428154
				0.8 1.53740665143427
				0.9 1.71670832559655
				1 1.84184998584447
				1.1 1.91848290198952
				1.2 1.95869248491735
				1.3 1.97612958981966
				1.4 1.98203303167729
				1.5 1.98345064795612
			};
			\addplot [semithick, goldenrod18818934]
			table {%
				0 0
				0.1 -0.142284418582873
				0.2 -0.327167755100247
				0.3 -0.508511917110884
				0.4 -0.653964989998738
				0.5 -0.733505955781393
				0.6 -0.723433547704958
				0.7 -0.612121157413337
				0.8 -0.404499179108363
				0.9 -0.123005629452394
				1 0.196219716808048
				1.1 0.511455421580314
				1.2 0.78557116233474
				1.3 0.994833070935357
				1.4 1.13313358827929
				1.5 1.21017056119002
			};
			\addplot [semithick, darkturquoise23190207]
			table {%
				0 0
				0.1 0.0609887782304123
				0.2 0.150907393869079
				0.3 0.257921151080575
				0.4 0.377676185548112
				0.5 0.506344158533972
				0.6 0.639670870884963
				0.7 0.772717021368166
				0.8 0.899858957875265
				0.9 1.01495403041967
				1 1.11163606958481
				1.1 1.18371461204166
				1.2 1.22564494376698
				1.3 1.23302556716269
				1.4 1.20307019361127
				1.5 1.13499602496936
			};
			\addplot [semithick, steelblue31119180]
			table {%
				0 0
				0.1 -1.22865806417849e-09
				0.2 -3.37370942293802e-08
				0.3 -3.51000746690033e-07
				0.4 -2.15121648727834e-06
				0.5 -9.4172499480738e-06
				0.6 -3.2624462286357e-05
				0.7 -9.50889264341725e-05
				0.8 -0.000242583519385775
				0.9 -0.000556528853548066
				1 -0.00117059065575531
				1.1 -0.00228993131249633
				1.2 -0.00421173754886717
				1.3 -0.00734505794461648
				1.4 -0.0122275159187804
				1.5 -0.0195361921019912
			};
			\addplot [semithick, darkorange25512714]
			table {%
				0 0
				0.1 0.000879139698502285
				0.2 0.00391276792186035
				0.3 0.011904137149789
				0.4 0.0292987857104662
				0.5 0.0621169867491682
				0.6 0.11739050950704
				0.7 0.202105802767826
				0.8 0.321789253373984
				0.9 0.478989473011206
				1 0.671978919743209
				1.1 0.893991835198929
				1.2 1.13323460742101
				1.3 1.37376384722078
				1.4 1.59715687670533
				1.5 1.78473646121355
			};
			\addplot [semithick, forestgreen4416044]
			table {%
				0 0
				0.1 -0.0556328237041062
				0.2 -0.0644071305415659
				0.3 -0.03935490306558
				0.4 0.0144982122074191
				0.5 0.0929386756650732
				0.6 0.191861791396003
				0.7 0.30714383374354
				0.8 0.434667307051692
				0.9 0.570397705485222
				1 0.71047577458536
				1.1 0.85130894488201
				1.2 0.98965320639634
				1.3 1.12268034459132
				1.4 1.24802760300312
				1.5 1.36382832245457
			};
			\addplot [semithick, crimson2143940]
			table {%
				0 0
				0.1 -0.141175458606859
				0.2 -0.306080538887106
				0.3 -0.478645665858656
				0.4 -0.648674951755348
				0.5 -0.805836684865538
				0.6 -0.93974402984178
				0.7 -1.04068210391287
				0.8 -1.10051620076344
				0.9 -1.1135753881019
				1 -1.07735464147189
				1.1 -0.992910692099669
				1.2 -0.864870512252902
				1.3 -0.701028435087886
				1.4 -0.511571345895449
				1.5 -0.308030714939454
			};
			\addplot [semithick, mediumpurple148103189]
			table {%
				0 0
				0.1 -0.000900153852729457
				0.2 -0.00870033169401274
				0.3 -0.0328236077346248
				0.4 -0.0828465798432786
				0.5 -0.166490292583876
				0.6 -0.288182937782573
				0.7 -0.448180092969698
				0.8 -0.64226294878014
				0.9 -0.861996568002043
				1 -1.09548136744616
				1.1 -1.32848742334454
				1.2 -1.54582979004324
				1.3 -1.73282716031885
				1.4 -1.87668687856215
				1.5 -1.96767561312314
			};
			\addplot [semithick, sienna1408675]
			table {%
				0 0
				0.1 0.0428039345849833
				0.2 0.200359138694995
				0.3 0.43116635378086
				0.4 0.70072633331012
				0.5 0.978100950949695
				0.6 1.23823459955363
				0.7 1.463710391287
				0.8 1.6452533722017
				0.9 1.78099755374911
				1 1.87482860483503
				1.1 1.93424667051794
				1.2 1.96820723262288
				1.3 1.98531869686807
				1.4 1.99263569174303
				1.5 1.995123671452
			};
			\addplot [semithick, orchid227119194]
			table {%
				0 0
				0.1 0.0035686952118551
				0.2 0.00905312771276279
				0.3 0.0169132091450601
				0.4 0.0278540127262042
				0.5 0.042717379224593
				0.6 0.0624616938333775
				0.7 0.088136655591845
				0.8 0.120846744570928
				0.9 0.161702072641713
				1 0.211757087072728
				1.1 0.271938973867927
				1.2 0.34296889764839
				1.3 0.425280395943974
				1.4 0.518940224247505
				1.5 0.623577616278161
			};
			\addplot [semithick, gray127]
			table {%
				0 0
				0.1 -0.0422206206891298
				0.2 -0.102932012705155
				0.3 -0.172488600555219
				0.4 -0.247361778718027
				0.5 -0.325026965967674
				0.6 -0.403256816843088
				0.7 -0.479903605561202
				0.8 -0.552826480097934
				0.9 -0.619880133156566
				1 -0.678936360794947
				1.1 -0.727926072404172
				1.2 -0.764894867087227
				1.3 -0.78806730497325
				1.4 -0.7959156686943
				1.5 -0.787229147505629
			};
			\addplot [semithick, goldenrod18818934]
			table {%
				0 0
				0.1 0.0241302679817232
				0.2 0.0680946167899097
				0.3 0.131614576981477
				0.4 0.213782846219197
				0.5 0.313078408696439
				0.6 0.427452060650156
				0.7 0.554404754629214
				0.8 0.691069855293843
				0.9 0.834302670606155
				1 0.980777395308958
				1.1 1.1270901440092
				1.2 1.26986602708308
				1.3 1.40586788509283
				1.4 1.53210423499217
				1.5 1.64593418786568
			};
			\addplot [semithick, darkturquoise23190207]
			table {%
				0 0
				0.1 0.0634961428409717
				0.2 0.150936056630778
				0.3 0.251115092221971
				0.4 0.359962705062319
				0.5 0.474297864202164
				0.6 0.591128734291673
				0.7 0.707517159764099
				0.8 0.820604413762578
				0.9 0.92770017127417
				1 1.02639757266592
				1.1 1.11469403474842
				1.2 1.1911030547403
				1.3 1.25474461457096
				1.4 1.30540337704016
				1.5 1.34354553220006
			};
			\addplot [semithick, steelblue31119180]
			table {%
				0 0
				0.1 0.0064209036987064
				0.2 0.0113264548156504
				0.3 0.0134791438343954
				0.4 0.0140005001677935
				0.5 0.014036946593645
				0.6 0.0140298159162089
				0.7 0.0136714085361197
				0.8 0.0108560422630213
				0.9 -8.94605103713975e-05
				1 -0.0294622107120939
				1.1 -0.0916750159807492
				1.2 -0.202666469246076
				1.3 -0.37533185111259
				1.4 -0.613883627223936
				1.5 -0.908779976361372
			};
			\addplot [semithick, darkorange25512714]
			table {%
				0 0
				0.1 -0.154770187563403
				0.2 -0.34059872132443
				0.3 -0.538536319413607
				0.4 -0.735519927840831
				0.5 -0.917860725069764
				0.6 -1.0716701056547
				0.7 -1.18404528986862
				0.8 -1.24444688321805
				0.9 -1.24598644466329
				1 -1.18640073177
				1.1 -1.0685364431813
				1.2 -0.9002354859334
				1.3 -0.693596102370547
				1.4 -0.463679769939109
				1.5 -0.226823455188256
			};
			\addplot [semithick, forestgreen4416044]
			table {%
				0 0
				0.1 0.00792587039031278
				0.2 0.0113280595545634
				0.3 0.0125081599321268
				0.4 0.0128180473657437
				0.5 0.0128716498630801
				0.6 0.0128762044751562
				0.7 0.0128762925077752
				0.8 0.0128762922817481
				0.9 0.0128760144294598
				1 0.0128669031443374
				1.1 0.0127783842724767
				1.2 0.0123151190544829
				1.3 0.0106479282994432
				1.4 0.00598930507522744
				1.5 -0.00483419851513367
			};
			\addplot [semithick, crimson2143940]
			table {%
				0 0
				0.1 0.109171804663082
				0.2 0.22451883344933
				0.3 0.325072613499141
				0.4 0.399780066065095
				0.5 0.440218114518341
				0.6 0.440425779092495
				0.7 0.397296782663336
				0.8 0.310923083266793
				0.9 0.184684339320266
				1 0.0249973717887941
				1.1 -0.159283063092149
				1.2 -0.357754152222951
				1.3 -0.559532857932588
				1.4 -0.754354279780215
				1.5 -0.933545048140873
			};
			\addplot [semithick, mediumpurple148103189]
			table {%
				0 0
				0.1 0.0635098802227673
				0.2 0.167121824683439
				0.3 0.319146096821256
				0.4 0.520399545316902
				0.5 0.763611700885051
				0.6 1.03335165839915
				0.7 1.30755200507466
				0.8 1.5604362815715
				0.9 1.76631676236122
				1 1.90353781448122
				1.1 1.95781252301026
				1.2 1.92433336651605
				1.3 1.80829530763252
				1.4 1.62379091076413
				1.5 1.39135104649475
			};
			\addplot [semithick, sienna1408675]
			table {%
				0 0
				0.1 0.0745811594404885
				0.2 0.1135456888571
				0.3 0.13961271583762
				0.4 0.15771305436885
				0.5 0.170401030159054
				0.6 0.179278386747667
				0.7 0.185438076106486
				0.8 0.189657243472115
				0.9 0.192499166698319
				1 0.194374659971803
				1.1 0.195582644972433
				1.2 0.196338747118438
				1.3 0.196796345285722
				1.4 0.197062490304757
				1.5 0.197210115630448
			};
			\addplot [semithick, orchid227119194]
			table {%
				0 0
				0.1 0.0868128833375501
				0.2 0.186299462177149
				0.3 0.282405780535739
				0.4 0.369479647534353
				0.5 0.444558498131856
				0.6 0.50639390275349
				0.7 0.555040606544013
				0.8 0.591529214804645
				0.9 0.617539864995043
				1 0.635082124273478
				1.1 0.646209607957742
				1.2 0.652797819949039
				1.3 0.656402704214319
				1.4 0.658202024684651
				1.5 0.6590070855651
			};
			\addplot [semithick, gray127]
			table {%
				0 0
				0.1 -0.0400836466355645
				0.2 -0.0277757459403751
				0.3 0.0704686152864151
				0.4 0.261153271305897
				0.5 0.529321483216533
				0.6 0.842392567288236
				0.7 1.15874585503001
				0.8 1.4392523164426
				0.9 1.65765982079686
				1 1.80580692478994
				1.1 1.8919462412206
				1.2 1.93376339830921
				1.3 1.95002516120292
				1.4 1.9547639090777
				1.5 1.9556869035219
			};
			\addplot [semithick, goldenrod18818934]
			table {%
				0 0
				0.1 1.64412987680995e-05
				0.2 0.000113296145348274
				0.3 0.000496038825228973
				0.4 0.00165792893539729
				0.5 0.00457700781156722
				0.6 0.0109261373883061
				0.7 0.023237623924067
				0.8 0.0449528871400761
				0.9 0.0802937438314276
				1 0.133916759566411
				1.1 0.21035254145555
				1.2 0.313280190526733
				1.3 0.444733031504438
				1.4 0.604364385757051
				1.5 0.788912822904061
			};
			\addplot [semithick, darkturquoise23190207]
			table {%
				0 0
				0.1 0.0642964554513866
				0.2 0.131902587226537
				0.3 0.208727789982506
				0.4 0.294857061758689
				0.5 0.389662394084107
				0.6 0.492139748787922
				0.7 0.601017128517636
				0.8 0.71481807630079
				0.9 0.831913536204539
				1 0.950570113665182
				1.1 1.06899665848479
				1.2 1.18538924016563
				1.3 1.29797397701209
				1.4 1.40504700752743
				1.5 1.50501089475875
			};
		\end{axis}
		
	\end{tikzpicture}

	\caption{The time dependence function $f(\cdot, K_j)$ for a few instances $K_j$ of the set of parameter distributions.}
	\label{fig:f_time_dep}
\end{figure}

For MRIs of the brain we can combine each image with synthetic rigid deformations. We choose time-dependent rotations and shifts according to the following parameters: 
\begin{equation}
	\bar{f}
	(t,a,b,p,q,C) = \mathrm{sign}(p)(2|\mathrm{sign}(p)\sin(a+Cb(t/T)^q)/2+0.5|^{1+|p|} - 1) \nonumber
\end{equation}
\begin{align}
	&\alpha \in \mathcal{U}(-\pi/36, \pi/36) , \qquad 
	d_\mathrm{x},d_\mathrm{y} \in \mathcal{U}(-3, 3) , \qquad
	q \in \mathcal{U}(1/1.3, 1.3) ,\nonumber \\
	&p \in \mathcal{U}(-4, 4), \qquad 
	a \in \mathcal{U}(0, 2\pi) , \qquad 
	b \in \mathcal{U}(0.5, 2), \qquad
	C = 1.5 \nonumber
\end{align}
\begin{equation}
	U_t[s] = T_t[R_t[s]], \qquad T_t[s](x,y) = s(x + f(t,K_1)d_\mathrm{x}, y + f(t, K_2)d_\mathrm{y}), \qquad R_t = R^{f(t, K_3)\alpha} \nonumber 
\end{equation}
with $f(t, K_j) = \bar{f}(t,K_j) - \bar{f}(0,K_j)$ where $K_j=\{a,b,p,q,C,d_\mathrm{x},d_\mathrm{y},\alpha \}$ is an instance of the set of parameter distributions, with the rotation matrix $R^\alpha$ of angle $\alpha$ and with the uniform distribution $\mathcal{U}$. 
A few instances of the time dependence function $f$ are plotted in figure \ref{fig:f_time_dep}.

All parameters are chosen probabilistically. Recall, that we model data acquisition to start at $t=0$ and end at $t=T$. 
The shifts $d_x$ and $d_y$ are to be understood as percent of the full field of view.

\subsubsection{Abdominal MRI}\label{sec:def_abdominal_MRI}
The simulation of breathing motion is much more complicated. Respiratory motion can be characterized as an elastic deformation of the body, modelled by the Navier-Cauchy equations\cite{hahn2022}. Therefore, on a set of 164 abdominal images, we generate a segmentation of the full abdominal cavity $\Omega$ by hand. Using these ground truths, we train a simple segmentation network (a U-Net\cite{ronneberger2015} with depth 4 and 16 channels) to extrapolate the segmentation to the rest of the data set. Using these segmentations, inside the cavity, we simulate the static Navier-Cauchy equations
\begin{equation}
	0 = \mu \Delta \vec{u} + (\lambda+\mu)\nabla(\nabla\cdot \vec{u}) + \rho \vec{k} \nonumber
\end{equation}
by applying a cranial-caudal volumetric force $\vec{k}$, assuming $\rho \equiv 1$. At the boundary we choose free-slip boundary conditions
\begin{equation}
	\vec{n} \cdot \vec{u} = 0, \quad \frac{\partial \vec{u}}{\partial \vec{n}} = 0 \qquad \forall \vec{u} \in \partial \Omega .
\end{equation}
The force is temporarily varied by the time dependence function $f(t,K_4)$ defined in \ref{sec:def_brain_MRI}, with an instance $K_4$ of the set of parameter distributions, and a scaling factor $c \in \mathcal{U}(0,1)$. Additionally, we added globally rigid deformations, as previously used for the case of brain MRIs, however weighted by a factor of $0.5$.

\section{Results}

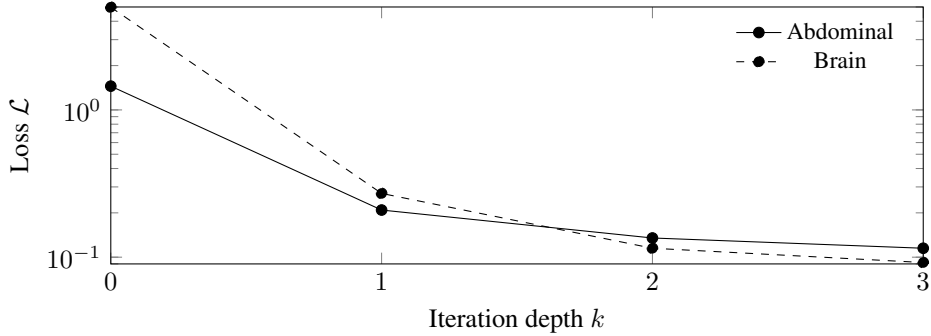
\begin{figure}
	\centering
	\begin{tikzpicture}
		
		\begin{semilogyaxis}[
			xmin = 0, xmax = 3, ymin = 0.09, ymax = 5,height=5cm,width = 0.75\textwidth, xlabel=Iteration depth $k$, ylabel=Loss $\mathcal{L}$, legend style={nodes={scale=0.9},draw=none,fill=none}, xtick={0,1,2,3}]
			\addplot[black, mark=*] table[x=iter, y=abdom] {B.dat};
			\addplot[black, mark=*, dashed] table[x=iter, y=brain] {B.dat};
			\legend{Abdominal, Brain}
		\end{semilogyaxis}
		
	\end{tikzpicture}
	\caption{Test Error Loss performance of proposed network over iteration depth. At iteration number zero the loss of an identity deformation is shown. }
	\label{fig:loss_performance}
\end{figure}

We evaluate the performance of the proposed network on simulated brain and abdominal HASTE MRIs. Therefore, we pair each static image with at least 5 different deformation fields. This significantly prevents overfitting, i.e. the network is not able to learn a direct relation from the image to a deformation, but rather has to assess the motion artefacts within the image.
For training, we split the dataset into 80\% training, 10\% validation and 10\% test data. All presented results and images are produced with the test data set. The computation time of the proposed method is compared to the Newton minimization scheme for rigid deformations as presented by Cordero and the GAN presented by Usman in table \ref{tab:computation_times}. With an effective computation time of $1.5$ seconds per slice, a clinical application is feasible. The loss evolution over cascades is plotted in figure \ref{fig:loss_performance} showing successive reduction. 

\subsection{Brain MRI}
\begin{figure}
	\centering
	\begin{subfigure}[b]{0.25\textwidth}
		\centering
		\includegraphics[width=\textwidth]{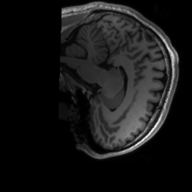}
		\caption*{ground truth $s^\mathrm{ref}$}
\end{subfigure}\hfil
\begin{subfigure}[b]{0.25\textwidth}
	\centering
	\includegraphics[width=\textwidth]{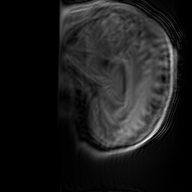}
	\caption*{static reconstruction}
\end{subfigure}\hfil
\begin{subfigure}[b]{0.25\textwidth}
	\centering
	\includegraphics[width=\textwidth]{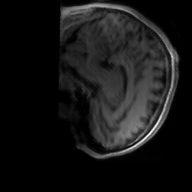}
	\caption*{GAN}
\end{subfigure}

\begin{subfigure}[b]{0.25\textwidth}
	\centering
	\includegraphics[width=\textwidth]{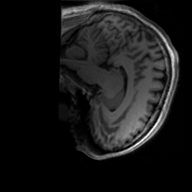}
	\caption*{proposed}
\end{subfigure}\hfil
\begin{subfigure}[b]{0.25\textwidth}
	\centering
	\includegraphics[width=\textwidth]{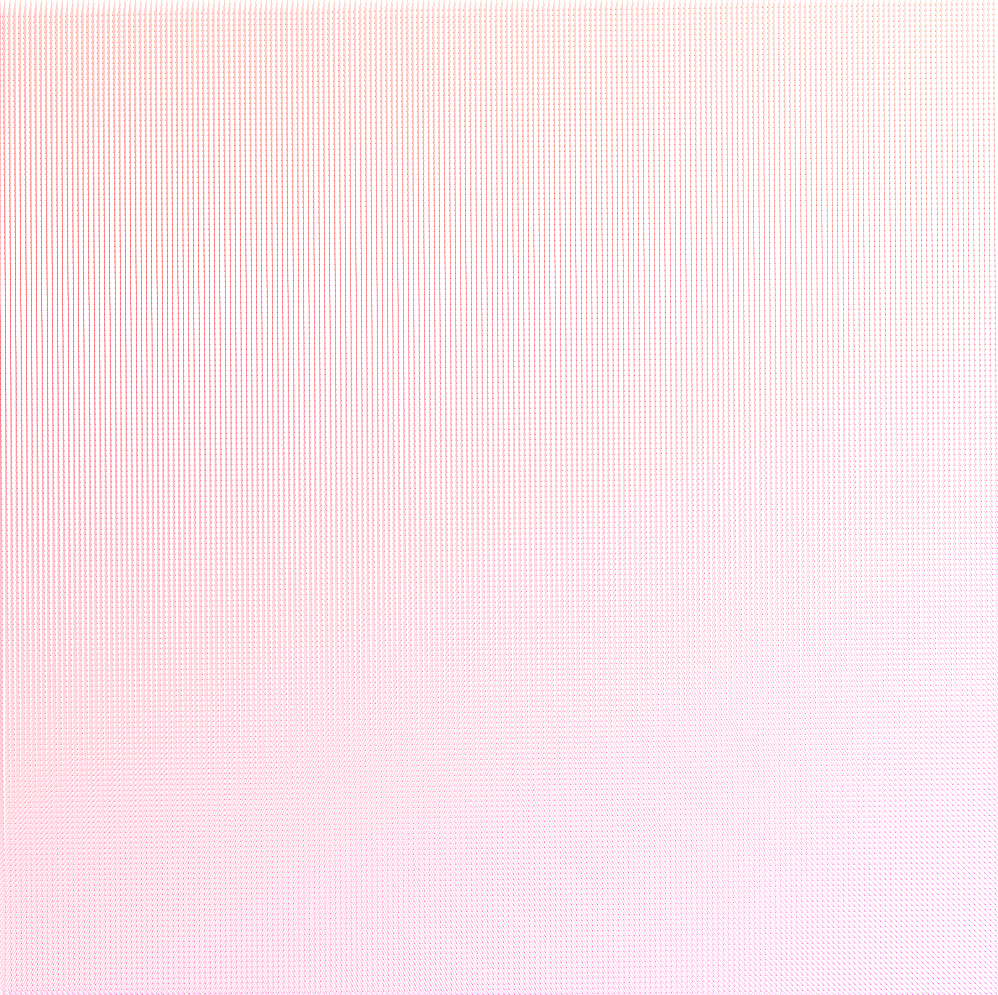}
	\caption*{ground truth $U_{7}$}
\end{subfigure}\hfil
\begin{subfigure}[b]{0.25\textwidth}
	\centering
	\includegraphics[width=\textwidth]{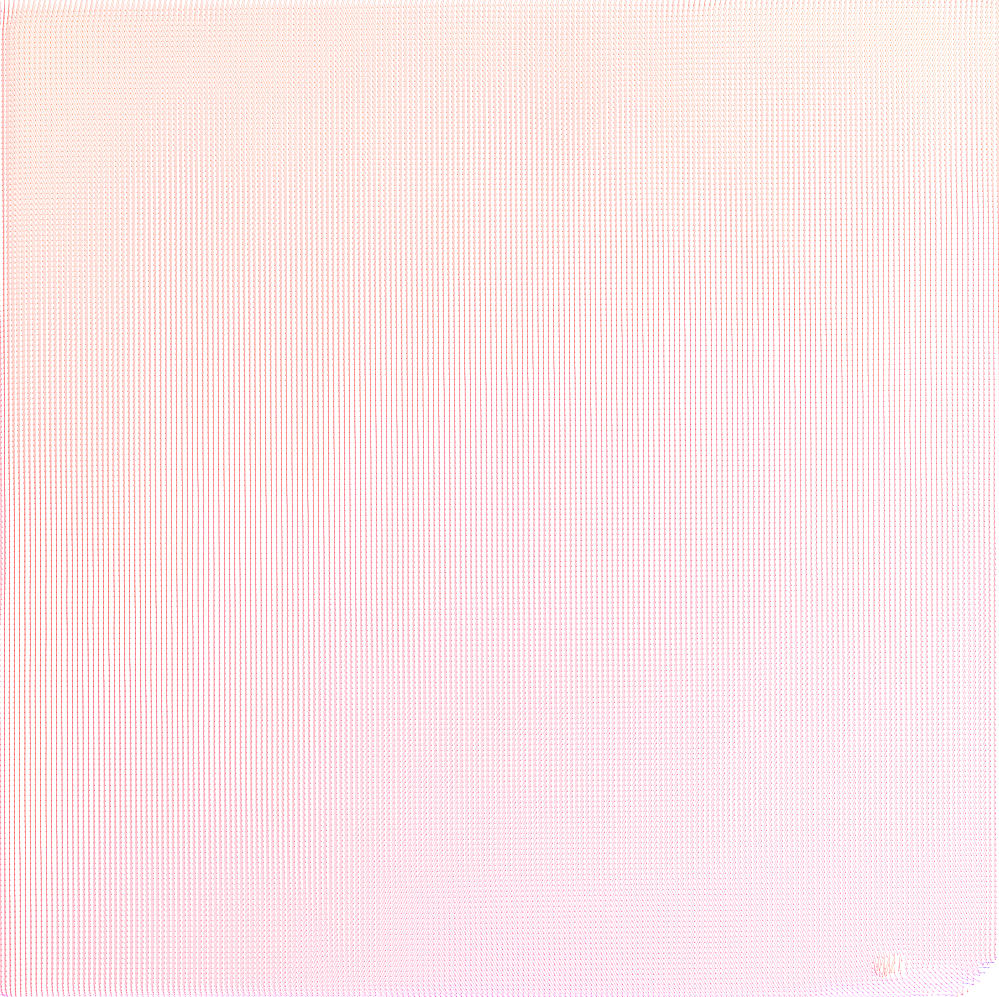}
	\caption*{prediction $U^\mathrm{est}_7$}
\end{subfigure}

\caption{Results on Brain HASTE MRIs $N^\mathrm{coils}=4$. All reconstructions have a resolution of $192\times192$. } 
\label{fig:reconstructions_synthetic_Brain}
\end{figure}

For MRIs of the brain we use 3600 images of the Yale dataset \cite{finn2015} and combine it with synthetic rigid deformations described in section \ref{sec:def_brain_MRI}. 
One instance of a reconstruction is shown in figure \ref{fig:reconstructions_synthetic_Brain}. The proposed method manages to reproduce the structure of the brain with much higher accuracy compared to the prediction of a GAN. The estimated deformation fields are accurate and well extrapolated to the exterior of the image. Further, the deformation fields are smooth and almost exactly rigid. Motion artefacts are successfully compensated and no limited data artefacts like streaks appear. Quantitative results for different noise levels are shown in the left column of table \ref{tab:quality_of_reconstructions}. The residual value remains almost unchanged, while the reconstruction quality of the proposed method significantly improves. Hence the residual value can indeed not be used for optimizations w.r.t. deformation fields. It can only serve as an indicator if data consistency is sufficiently preserved. The proposed method achieves the best scores in all error metrics (peak signal to noise ratio(PSNR), structural similarity index(SSIM) and mean squared error(MSE)) by a substantial margin.

\subsection{Abdominal MRI}
\begin{figure}
\centering
\begin{subfigure}[b]{0.25\textwidth}
	\centering
	\includegraphics[width=\textwidth]{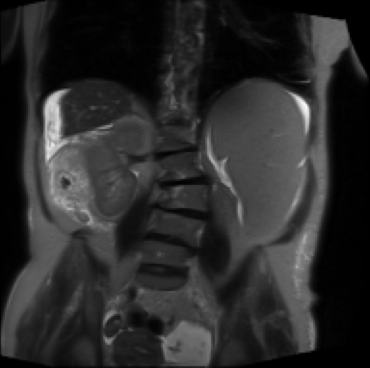}
	\caption*{ground truth $s^\mathrm{ref}$}
\end{subfigure}\hfil
\begin{subfigure}[b]{0.25\textwidth}
\centering
\includegraphics[width=\textwidth]{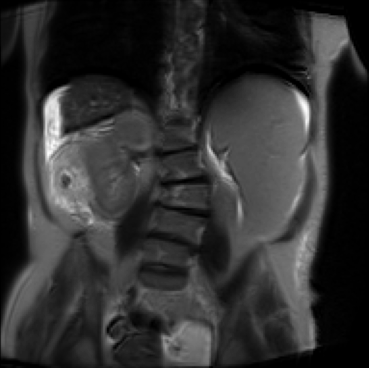}
\caption*{static reconstruction}
\end{subfigure}\hfil
\begin{subfigure}[b]{0.25\textwidth}
\centering
\includegraphics[width=\textwidth]{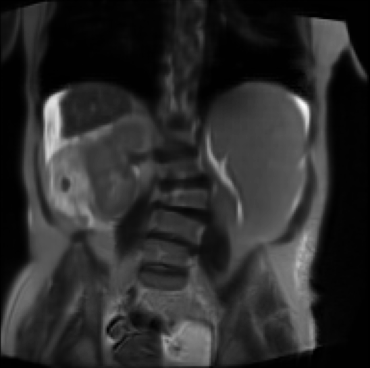}
\caption*{GAN}
\end{subfigure}

\begin{subfigure}[b]{0.25\textwidth}
\centering
\includegraphics[width=\textwidth]{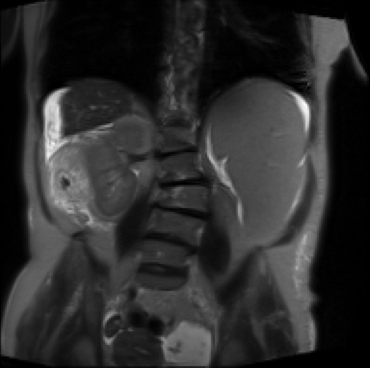}
\caption*{proposed}
\end{subfigure}\hfil
\begin{subfigure}[b]{0.25\textwidth}
\centering
\includegraphics[width=\textwidth]{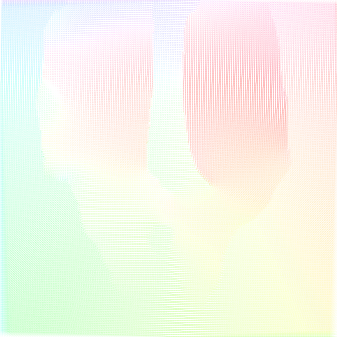}
\caption*{ground truth $U_{7}$}
\end{subfigure}\hfil
\begin{subfigure}[b]{0.25\textwidth}
\centering
\includegraphics[width=\textwidth]{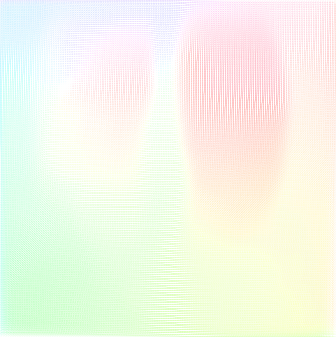}
\caption*{prediction $U^\mathrm{est}_7$}
\end{subfigure}

\caption{Results on Abdominal HASTE MRIs $N^\mathrm{coils}=4$. All reconstructions have a resolution of $192\times192$. } 
\label{fig:reconstructions_synthetic_Abdominal}
\end{figure}

We are equipped with a data set of 2400 liver MRIs showing the full abdomen of a patient. We pair these images with simulated breathing deformations computed by solving the static Navier-Cauchy equations described in section \ref{sec:def_abdominal_MRI}.
An instance of the resulting deformation fields can be seen in figure \ref{fig:reconstructions_synthetic_Abdominal}. The simulated deformation fields are visually similar to what 3D registration algorithms produce c.f\cite{huttinga2022}. The estimated deformation is smooth at the discontinuities of the ground truth. This can be interpreted as an uncertain determination of the edges of the abdominal cavity. More accurate segmentations can most likely lead to improved results. The overall appearance is however excellent and breathing motion is detected reliably. Motion artefacts are well compensated for and details are recovered. The standard GAN reconstruction manages to remove dominant motion artefacts but removes details at the same time. Quantitative results on different noise levels are shown in the right column of table \ref{tab:quality_of_reconstructions}. The proposed method achieves the best scores in all error metrics.

\begin{table}
	\centering
	\begin{tabular}{ l||l|l|l } 
		\hline
		\hline
		
		Method & proposed & Cordero & Usman \\
		\hline
		
		\hline
		
		time         & $3.1$ & $93.8$ & $0.046$ \\
		
		\hline
		\hline
	\end{tabular}
	\caption{Computation time in seconds for one data sample and $N^\mathrm{coils}=4$. Usually, more than one data sample is acquired. This can benefit computation times due to parallelizability. A batch of 4 data samples requires $6.06s$ to be processed resulting in $1.5s$ per data sample for the proposed method. We used a workstation with Intel i9-10900X, 128GB Ram and NVIDIA RTX 3090. }
	\label{tab:computation_times}
\end{table}

\begin{table}
	\centering
	\begin{tabular}{ l||cccc|cccc } 
		\hline
		\hline
		
		Method & Res & PSNR & SSIM & MSE & Res & PSNR & SSIM & MSE \\
		\hline
		& \multicolumn{4}{c}{Brain}                                & \multicolumn{4}{c}{Abdominal} \\
		& \multicolumn{8}{c}{noise level: 0\%} \\
		\hline
		
		static       & 9.330 & 30.66 & 0.899 & 6.746                                      & \textbf{9.378} & 33.16 & 0.908 & 4.706\\
		GAN (Usman)  & 33.83 & 32.00 & 0.938 & 5.621                                      & 28.34 & 32.61 & 0.932 & 4.647\\
		CG($y, U^\mathrm{est}$)     & \textbf{8.108} & 39.73 & 0.987 & 2.191	& 10.56 &  36.45 & 0.957 & 3.081\\
		CG($y, U^\mathrm{est}$) + additional GAN     & 14.92 & \textbf{41.67} & \textbf{0.992} & \textbf{1.710}  & 20.05 & \textbf{ 36.62} & \textbf{0.968} & \textbf{2.933}\\
		\hline

		& \multicolumn{8}{c}{noise level: 5\%} \\
		\hline

		static       & 21.81 & 30.77 & 0.892 & 6.858                                      & \textbf{26.68} & 32.22 & 0.887 & 5.116\\
		GAN (Usman)  & 36.17 & 32.29 & 0.935 & 5.555                                      & 35.89 & 32.32 & 0.926 & 4.800\\
		CG($y, U^\mathrm{est}$)      & \textbf{21.79} & 33.62 & 0.934 & 4.453  & 26.76 & 34.41 & 0.917 & 3.745\\
		CG($y, U^\mathrm{est}$) + additional GAN      & 29.23 & \textbf{37.42} & \textbf{0.972} & \textbf{2.853}  & 30.29 & \textbf{35.89} & \textbf{0.954} & \textbf{3.133}\\
		\hline

		\hline
	\end{tabular}
	\caption{Quantitative comparison of reconstruction quality. Here, CG($y, U^\mathrm{est}$) is the CG-SENSE method with data $y$ and estimated deformation fields $U^\mathrm{est}$ by the proposed method. Remaining artefacts can be compensated by the application of an additional GAN.}         
	\label{tab:quality_of_reconstructions}
\end{table}


\section{Conclusion}

We presented a framework to compensate 2D in-plane motion artefacts in HASTE MRIs. The proposed method does not require reference scans but can predict the intra-image deformations from inconsistencies alone.

The results surpass standard GAN approaches of the literature while maintaining data consistency. Even though the proposed method is specifically designed to cope with MRI sequences without redundancies, it can still be applied to sequences with redundancies like CINE.

This supervised approach requires suitably many known deformation fields. Unsupervised training w.r.t. $U$  can be achieved by the implementation of the loss function

\begin{equation}
\mathcal{L}^\mathrm{unsup}(U, s^\mathrm{ref}) = |s^\mathrm{ref} - \mathcal{R}(y,U)|_2^2 + |\nabla U|_2^2,
\end{equation}

where $\mathcal{R}$ is the dynamic reconstruction operator. Note, that this approach requires gradients of an iterative reconstruction operator. Further, the loss function is reference state dependent, which can hinder the convergence of the training and has to be addressed in future work.

The proposed strategy can also be applied to other imaging modalities like computerized tomography (CT) as well as to 3D applications. In case of semi-3D acquisitions, i.e. multiple 2D slices, the extension is straight forward. For the case of actual 3D Fourier coefficients, a more sophisticated adaption to this scenario is required. First, the GAN prediction has to be a 3D model. Using this prediction, we can still compute partial reconstructions $s_i^\mathrm{GAN}$. Since the corresponding data $y_i^\mathrm{GAN}$ to $s_i^\mathrm{GAN}$ is only a stripe of a slice in 3D-k-space, these data do not provide information in the third dimension. Therefore, we can reduce its dimension to 2D. The same holds for $s_i$, leading to a 2D to 2D registration task with the goal of estimating 3D deformation fields. Clearly, the degrees of freedom limit the complexity of these deformations significantly. The complex argument allows to detect through-plane motion by means of rigid deformations. Since complex network structures are required for this task, results are subject to future research. Motion of further complexity requires more information in the third dimension limiting temporal resolution of deformation estimates or even requires sequence design. Due to results in the literature for slice-to-volume registration\cite{ferrante2017}, we expect that one can successfully extend the proposed idea to 3D whole brain MRIs acquired with cartesian sequences using partial Fourier and parallel imaging acceleration.

\bibliographystyle{unsrt}
\bibliography{Literatur}

\begin{thebibliography}{10}

\bibitem{cordero2016}
L.~{Cordero-Grande}, R.~P. A.~G. {Teixeira}, E.~J. {Hughes}, J.~{Hutter}, A.~N.
  {Price}, and J.~V. {Hajnal}.
\newblock Sensitivity encoding for aligned multishot magnetic resonance
  reconstruction.
\newblock {\em IEEE Trans. Comput. Imaging}, 2(3):266--280, 2016.

\bibitem{burger2018}
Martin Burger, Hendrik Dirks, and Carola-Bibiane Sch\"onlieb.
\newblock A variational model for joint motion estimation and image
  reconstruction.
\newblock {\em SIAM J. Imaging Sci.}, 11(1):94--128, 2018.

\bibitem{feinler2023}
Mathias~S. Feinler and Bernadette~N. Hahn.
\newblock Retrospective motion correction in gradient echo mri by explicit
  motion estimation using deep cnns, 2023.

\bibitem{andre2015}
Jalal~B. Andre, Brian~W. Bresnahan, Mahmud Mossa-Basha, Michael~N. Hoff,
  C.~Patrick Smith, Yoshimi Anzai, and Wendy~A. Cohen.
\newblock Toward quantifying the prevalence, severity, and cost associated with
  patient motion during clinical mr examinations.
\newblock {\em Journal of the American College of Radiology}, 12(7):689--695,
  2015.

\bibitem{pruessmann1999}
Klaas~P Pruessmann, Markus Weiger, Markus~B Scheidegger, and Peter Boesiger.
\newblock Sense: sensitivity encoding for fast mri.
\newblock {\em Magn. Reson. Med.}, 42(5):952--962, 1999.

\bibitem{allison2012}
Michael~J. Allison, Sathish Ramani, and Jeffrey~A. Fessler.
\newblock Accelerated regularized estimation of mr coil sensitivities using
  augmented lagrangian methods.
\newblock {\em IEEE Trans. Med. Imaging}, 32(3):556--564, 2013.

\bibitem{uecker2014}
Martin Uecker, Peng Lai, Mark~J. Murphy, Patrick Virtue, Michael Elad, John~M.
  Pauly, Shreyas~S. Vasanawala, and Michael Lustig.
\newblock Espirit—an eigenvalue approach to autocalibrating parallel mri:
  Where sense meets grappa.
\newblock {\em Magn. Reson. Med.}, 71(3):990--1001, 2014.

\bibitem{ma2015}
Ya-Jun Ma, Wentao Liu, Xin Tang, and Jia-Hong Gao.
\newblock Improved sense imaging using accurate coil sensitivity maps generated
  by a global magnitude-phase fitting method.
\newblock {\em Magnetic Resonance in Medicine}, 74(1):217--224, 2015.

\bibitem{patel1997}
M.~R. Patel, R.~A. Klufas, R.~A. Alberico, and R.~R. Edelman.
\newblock Half-fourier acquisition single-shot turbo spin-echo (haste) mr:
  comparison with fast spin-echo mr in diseases of the brain.
\newblock {\em AJNR. American journal of neuroradiology}, 18(9):1635--40, 1997.

\bibitem{usman2020}
Muhammad Usman, Siddique Latif, Muhammad Asim, Byoung-Dai Lee, and Junaid
  Qadir.
\newblock Retrospective motion correction in multishot mri using generative
  adversarial network.
\newblock {\em Scientific Reports}, 10:4786, 03 2020.

\bibitem{armanious2018b}
Karim Armanious, Sergios Gatidis, Konstantin Nikolaou, Bin Yang, and Thomas
  Kustner.
\newblock Retrospective correction of rigid and non-rigid mr motion artifacts
  using gans.
\newblock In {\em IEEE ISBI (2019)}, pages 1550--1554, 2019.

\bibitem{goodfellow2016}
Ian Goodfellow, Yoshua Bengio, and Aaron Courville.
\newblock {\em Deep Learning}.
\newblock MIT Press, 2016.
\newblock \url{http://www.deeplearningbook.org}.

\bibitem{yoshida2022}
Nobukiyo Yoshida, Hajime Kageyama, Hiroyuki Akai, Koichiro Yasaka, Haruto
  Sugawara, Yukinori Okada, and Akira Kunimatsu.
\newblock Motion correction in mr image for analysis of vsrad using generative
  adversarial network.
\newblock {\em PLOS ONE}, 17(9):1--18, 09 2022.

\bibitem{bao2022}
Qinjia Bao, Yalei Chen, Chongxin Bai, Pingan Li, Kewen Liu, Zhao Li, Zhi Zhang,
  Jie Wang, and Chaoyang Liu.
\newblock Retrospective motion correction for preclinical/clinical magnetic
  resonance imaging based on a conditional generative adversarial network with
  entropy loss.
\newblock {\em NMR in Biomedicine}, 35(12):e4809, 2022.

\bibitem{johnson2019}
Patricia~M. Johnson and Maria Drangova.
\newblock Conditional generative adversarial network for 3d rigid-body motion
  correction in mri.
\newblock {\em Magnetic Resonance in Medicine}, 82(3):901--910, 2019.

\bibitem{haskell2019}
Melissa~W. Haskell, Stephen~F. Cauley, Berkin Bilgic, Julian Hossbach,
  Daniel~N. Splitthoff, Josef Pfeuffer, Kawin Setsompop, and Lawrence~L. Wald.
\newblock Network accelerated motion estimation and reduction (namer):
  Convolutional neural network guided retrospective motion correction using a
  separable motion model.
\newblock {\em Magnetic Resonance in Medicine}, 82(4):1452--1461, 2019.

\bibitem{adler2017}
Jonas Adler and Ozan \"Oktem.
\newblock Solving ill-posed inverse problems using iterative deep neural
  networks.
\newblock {\em Inverse Probl.}, 33(12):124007, Nov 2017.

\bibitem{adler2018}
Jonas Adler and Ozan \"Oktem.
\newblock Learned primal-dual reconstruction.
\newblock {\em IEEE Trans. Med. Imaging}, 37(6):1322--1332, Jun 2018.

\bibitem{hammernik2018}
Kerstin Hammernik, Teresa Klatzer, Erich Kobler, Michael~P. Recht, Daniel~K.
  Sodickson, Thomas Pock, and Florian Knoll.
\newblock Learning a variational network for reconstruction of accelerated mri
  data.
\newblock {\em Magnetic Resonance in Medicine}, 79(6):3055--3071, 2018.

\bibitem{schlemper2017}
Jo~Schlemper, Jose Caballero, Joseph~V. Hajnal, Anthony Price, and Daniel
  Rueckert.
\newblock A deep cascade of convolutional neural networks for dynamic mr image
  reconstruction, 2017.

\bibitem{balakrishnan2019}
Guha Balakrishnan, Amy Zhao, Mert~R. Sabuncu, John Guttag, and Adrian~V. Dalca.
\newblock Voxelmorph: A learning framework for deformable medical image
  registration.
\newblock {\em IEEE Transactions on Medical Imaging}, 38(8):1788--1800, 2019.

\bibitem{munoz2022}
Camila Munoz, Haikun Qi, Gastao Cruz, Thomas Küstner, René~M. Botnar, and
  Claudia Prieto.
\newblock Self-supervised learning-based diffeomorphic non-rigid motion
  estimation for fast motion-compensated coronary mr angiography.
\newblock {\em Magnetic Resonance Imaging}, 85:10--18, 2022.

\bibitem{seegoolam2019}
Gavin Seegoolam, Jo~Schlemper, Chen Qin, Anthony Price, Jo~Hajnal, and Daniel
  Rueckert.
\newblock Exploiting motion for deep learning reconstruction of
  extremely-undersampled dynamic mri.
\newblock In {\em Medical Image Computing and Computer Assisted Intervention --
  MICCAI 2019}, pages 704--712, Cham, 2019. Springer International Publishing.

\bibitem{teed2020}
Zachary Teed and Jia Deng.
\newblock Raft: Recurrent all-pairs field transforms for optical flow.
\newblock In Andrea Vedaldi, Horst Bischof, Thomas Brox, and Jan-Michael Frahm,
  editors, {\em Computer Vision -- ECCV 2020}, pages 402--419, Cham, 2020.
  Springer International Publishing.

\bibitem{hammernik2021}
Kerstin Hammernik, Jiazhen Pan, Daniel Rueckert, and Thomas Küstner.
\newblock Motion-guided physics-based learning for cardiac mri reconstruction.
\newblock In {\em 2021 55th Asilomar Conference on Signals, Systems, and
  Computers}, pages 900--907, 2021.

\bibitem{pan2022}
Jiazhen Pan, Daniel Rueckert, Thomas Küstner, and Kerstin Hammernik.
\newblock Learning-based and unrolled motion-compensated reconstruction for
  cardiac mr cine imaging, 2022.

\bibitem{batchelor2006}
P.~G. Batchelor, D.~Atkinson, P.~Irarrazaval, D.~L.~G. Hill, J.~Hajnal, and
  D.~Larkman.
\newblock Matrix description of general motion correction applied to multishot
  images.
\newblock {\em Magn. Reson. Med.}, 54(5):1273--1280, 2005.

\bibitem{scherzer2008}
Otmar Scherzer, Markus Grasmair, Harald Grossauer, Markus Haltmeier, and Frank
  Lenzen.
\newblock Variational methods in imaging.
\newblock In {\em Applied Mathematical Sciences}, 2008.

\bibitem{horn1981}
Berthold~K.P. Horn and Brian~G. Schunck.
\newblock Determining optical flow.
\newblock {\em Int. J. Artif. Intell.}, 17(1):185--203, 1981.

\bibitem{ilg2016}
Eddy Ilg, Nikolaus Mayer, Tonmoy Saikia, Margret Keuper, Alexey Dosovitskiy,
  and Thomas Brox.
\newblock Flownet 2.0: Evolution of optical flow estimation with deep networks.
\newblock {\em CoRR}, 2016.

\bibitem{vaillant2014}
Ghislain Vaillant, Claudia Prieto, Christoph Kolbitsch, Graeme Penney, and
  Tobias Schaeffter.
\newblock Retrospective rigid motion correction in k-space for segmented radial
  mri.
\newblock {\em IEEE Transactions on Medical Imaging}, 33(1):1--10, 2014.

\bibitem{unser1995}
M.~Unser, P.~Thevenaz, and L.~Yaroslavsky.
\newblock Convolution-based interpolation for fast, high-quality rotation of
  images.
\newblock {\em IEEE Transactions on Image Processing}, 4(10):1371--1381, 1995.

\bibitem{yang2017}
Chao Yang, Xin Lu, Zhe Lin, Eli Shechtman, Oliver Wang, and Hao Li.
\newblock High-resolution image inpainting using multi-scale neural patch
  synthesis.
\newblock In {\em 2017 IEEE CVPR}, pages 4076--4084, 2017.

\bibitem{pruessmann2001}
Klaas~P. Pruessmann, Markus Weiger, Peter Börnert, and Peter Boesiger.
\newblock Advances in sensitivity encoding with arbitrary k-space trajectories.
\newblock {\em Magnetic Resonance in Medicine}, 46(4):638--651, 2001.

\bibitem{marquez2018}
Enrique~S. Marquez, Jonathon~S. Hare, and Mahesan Niranjan.
\newblock Deep cascade learning.
\newblock {\em IEEE Trans. Neur. Net. Lear.}, 29(11):5475--5485, 2018.

\bibitem{hahn2022}
Bernadette~N. Hahn, Melina-Loren~Kienle Garrido, Christian Klingenberg, and
  Sandra Warnecke.
\newblock Using the navier-cauchy equation for motion estimation in dynamic
  imaging.
\newblock {\em Inverse Problems and Imaging}, 16(5):1179--1198, 2022.

\bibitem{ronneberger2015}
Olaf Ronneberger, Philipp Fischer, and Thomas Brox.
\newblock U-net: Convolutional networks for biomedical image segmentation.
\newblock volume 9351, pages 234--241, 10 2015.

\bibitem{finn2015}
Emily~S. Finn, Xilin Shen, Dustin Scheinost, Jessica Huang, Marvin~M. Chun,
  Xenophon Papademetris, and Todd~R. Constable.
\newblock Functional connectome fingerprinting: identifying individuals using
  patterns of brain connectivity.
\newblock {\em Nat. Neurosci.}, 18(11):1664 -- 1671, 2015.

\bibitem{huttinga2022}
Niek R.~F. Huttinga, Tom Bruijnen, Cornelis A.~T. Van Den~Berg, and Alessandro
  Sbrizzi.
\newblock Real-time non-rigid 3d respiratory motion estimation for mr-guided
  radiotherapy using mr-motus.
\newblock {\em IEEE Transactions on Medical Imaging}, 41(2):332--346, 2022.

\bibitem{ferrante2017}
Enzo Ferrante and Nikos Paragios.
\newblock Slice-to-volume medical image registration: A survey.
\newblock {\em Medical Image Analysis}, 39:101--123, 2017.

\end{thebibliography}

\end{document}